\input amstex
\magnification=\magstep1 
\baselineskip=13pt
\documentstyle{amsppt}
\vsize=8.7truein
\CenteredTagsOnSplits \NoRunningHeads
\def\vl{\operatorname{vol}}
\def\per{\operatorname{per}}
\def\AA{\Cal A}
\def\HHH{\Cal H}

\def\RR{\Cal R}
\def\CC{\Cal C}
\def\xx{\bold x}
\def\yy{\bold y}
\def\ttt{\bold t}
\def\sss{\bold s}
\def\PP{\Cal P}
\def\HH{\bold H}

\def\Pr{\bold{Pr\thinspace }}

\topmatter
\title Asymptotic estimates for the number of contingency tables, integer flows,  and volumes of 
transportation polytopes \endtitle
\author Alexander Barvinok \endauthor
\address Department of Mathematics, University of Michigan, Ann Arbor,
MI 48109-1043, USA \endaddress
\email barvinok$\@$umich.edu \endemail
\thanks This research was partially supported by NSF Grant DMS 0400617. \endthanks
 \abstract We prove an asymptotic estimate for the number of $m \times n$ non-negative integer 
 matrices (contingency tables) with prescribed row 
  and column sums  and, more generally, for 
 the number of integer feasible flows in a network.  Similarly,
 we estimate the volume of the polytope of $m \times n$ non-negative real matrices
 with prescribed row and column sums. 
 Our estimates are solutions of convex optimization problems and hence can be computed efficiently.
 As a corollary, we show that  if  row sums $R=\left(r_1, \ldots, r_m \right)$ and column sums
 $C=\left(c_1, \ldots, c_n \right)$ with $r_1 + \ldots + r_m =c_1 + \ldots +c_n =N$ are sufficiently
 far from constant vectors,  then, asymptotically, in the uniform probability space of the $m \times n$
 non-negative integer matrices with the total sum $N$ of entries,
  the event consisting of the matrices with row sums $R$ and the event consisting 
 of the matrices with column sums $C$ are positively correlated.
  \endabstract
\date August  2008
\enddate
\endtopmatter

\document

\head 1. Introduction and main results \endhead

Let $m>1$ and $n>1$ be integers and let $R=\left(r_1, \ldots, r_m\right)$ and 
$C=\left(c_1, \ldots, c_n \right)$ be positive integer vectors such that 
$$\sum_{i=1}^m r_i =\sum_{j=1}^n c_j =N.$$
We are interested in the number $\#(R, C)$ of $m \times n$ non-negative integer matrices,
also known as {\it contingency tables},
with row sums $R$ and column sums $C$, called {\it margins}.
Computing or estimating numbers $\#(R, C)$ has attracted a lot of attention, because 
of the relevance of these numbers in statistics, see \cite{Goo76}, \cite{DE85}, 
combinatorics, representation theory, and elsewhere, see  \cite{DG85}, \cite{DG04}.
Of interest are asymptotic formulas, see \cite{BBK72}, \cite{Ben74} 
and most recent \cite{CM07a}, \cite{GM07},
algorithms with rigorous estimates of the performance guarantees, see \cite{DKM97}, \cite{Mor02}, 
\cite{CD03}, \cite{BLV04}, and heuristic approaches which may lack  formal justification 
but tend to work well in practice \cite{Goo76}, \cite{DE85},  \cite{C+05}.

Our first main result is as follows.

\proclaim{(1.1) Theorem} Let $R=\left(r_1, \ldots, r_m \right)$ and 
$C=\left(c_1, \ldots, c_n \right)$ be positive integer vectors such that 
$r_1 + \ldots + r_m = c_1 + \ldots + c_n =N$.
Let us define a function
$$\split F(&\xx, \yy)=\left(\prod_{i=1}^m x_i^{-r_i} \right) \left( \prod_{j=1}^n y_j^{-c_j} \right)
\left( \prod_{ij} {1 \over 1 - x_i y_j} \right) \\ &\quad \text{for} \quad
\xx=\left(x_1, \ldots, x_m \right) \quad \text{and} \quad \yy=\left(y_1, \ldots, y_n \right).
\endsplit$$
Then $F(\xx, \yy)$ attains its minimum
$$\rho=\rho(R,C)=\min  \Sb 0 < x_1, \ldots, x_m <1 \\ 0< y_1, \ldots, y_n <1 \endSb F(\xx, \yy)$$
on the open cube $0<x_i, y_j<1$ and 
for the number $\#(R,C)$ of non-negative integer $m \times n$ matrices with row 
sums $R$ and column sums $C$ we have
$$\rho \ \geq \ \#(R,C) \ \geq \ N^{-\gamma(m+n)} \rho,$$
where $\gamma >0$ is an absolute constant.
\endproclaim 

More precisely, the lower bound we prove is
$$\split \#(R,C) \ \geq \  \ &{\Gamma\left({m+n \over 2} \right) \over 2e^5  
\pi^{{m+n-2 \over 2}} mn (N+mn)} \left({2\over (mn)^2 (N+1) (N+mn)}\right)^{m+n-1} \\
&\quad \times \left(\prod_{i=1}^m {r_i^{r_i} \over r_i!} \right) \left( \prod_{j=1}^n {c_j^{c_j} \over c_j!} \right) 
{N! (N+mn)! (mn)^{mn} \over  N^N  (N+mn)^{N+mn}(mn)!} \rho(R,C)\endsplit$$
provided $m +n \geq 10$.
Recall that from Stirling's formula
$${s!  \over s^s} =\sqrt{2\pi s} \left(1 +O(s^{-1}) \right)$$
and hence the product in front of $\rho(R,C)$ indeed exceeds $N^{-\gamma(m+n)}$ for some 
absolute constant $\gamma>0$.

We note that the substitution $x_i=e^{-t_i}$ and $y_j=e^{-s_j}$ transforms the problem 
of computing $\rho$ into the problem of minimizing the convex function
$$\phi(\ttt,\sss)=\phi_{R, C}(\ttt,\sss)
=\sum_{i=1}^m r_i t_i + \sum_{j=1}^n c_j s_j -\sum_{ij} \ln\left(1-e^{-t_i -s_j} \right)$$
on the positive orthant $s_i, t_j>0$, so that methods of convex optimization can be applied 
to compute $\rho$ in time polynomial in $m+n$ and $\ln N$, see \cite{NN94}. 

Theorem 1.1 estimates the number $\#(R, C)$ of contingency tables within an $N^{O(m+n)}$
factor. This estimate provides, asymptotically,  the main term \break of  $\log \#(R,C)$ for all but very sparse cases,
 where margins $r_i$ and $c_j$ 
are small compared to the sizes $m$ and $n$ of the matrix. For example, if the margins 
$r_i$ and $c_j$ are at least linear in $m$ and $n$ then $\#(R, C)$ is at least as 
big as $\gamma^{mn}$ for some constant $\gamma>1$. By now, the sparse case of small 
$r_i$ and $c_j$ is well understood, thanks especially to the recent paper \cite{GM07}. 
The case of moderate to high margins seems to be the most difficult. To the author's knowledge, the 
estimate of Theorem 1.1 is the {\it only} rigorously proven effective estimate of $\#(R, C)$ for 
generic $R$ and $C$ (if all $r_i$'s are equal and all $c_j$'s are equal, recent paper 
\cite{CM07a} provides a precise asymptotic formula for the number of tables).
Theorem 1.1 allows us to find faults with the very intuitive ``independence heuristic'' for 
counting contingency tables 
and points out  at some strange ``attraction'' phenomena in the space of matrices. Quite 
counter-intuitively, we conclude that in the uniform probability space of the $m \times n$ non-negative 
integer matrices with the total sum of entries equal to $N$, the event consisting of the 
matrices with row sums $R$ and the event consisting of the matrices with column sums 
$C$ attract exponentially in $mn$ provided the vectors $R$ and $C$ are 
sufficiently far from constant vectors, see Section 2 for the precise statements and details.

Let us identify the space of $m \times n$ real matrices $X=\left(x_{ij}\right)$ with 
Euclidean space ${\Bbb R}^d$ for $d=mn$. In ${\Bbb R}^d$ we consider the {\it transportation 
polytope} $\PP=\PP(R, C)$ defined by the equations
$$\sum_{j=1}^n x_{ij} =r_i \quad \text{for} \quad i=1, \ldots, m, \qquad 
\sum_{i=1}^m x_{ij}=c_j \quad \text{for} \quad j=1, \ldots, n$$
and inequalities 
$$x_{ij} \geq 0 \quad \text{for all} \quad i,j.$$
As is known, $\PP$ is a polytope of dimension $(m-1)(n-1)$. We prove the 
following estimate for the volume of $\PP$, computed with respect to the Euclidean structure 
in the affine span of $\PP$, induced from ${\Bbb R}^d$. 

\proclaim{(1.2) Theorem} Let $R=\left(r_1, \ldots, r_m \right)$ and 
$C=\left(c_1, \ldots, c_n \right)$ be positive integer vectors such that 
$r_1 + \ldots + r_m = c_1 + \ldots + c_n =N$
and let $\PP =\PP(R,C)$ be the polytope of non-negative $m \times n$ matrices with 
row sums $r_1, \ldots, r_m$ and column sums $c_1, \ldots, c_n$.  

Let 
$$\beta=\beta(R,C)=\max \Sb X=\left(x_{ij}\right) \\ X \in \PP \endSb \prod_{ij} x_{ij}$$
be the maximum value of the product of entries of a matrix from $\PP$.
Then for the volume of $\PP$ we have
$$\beta e^{mn} N^{\gamma(m+n)}\  \geq \  \vl(\PP) \ \geq \   
\beta e^{mn} N^{-\gamma(m+n)} ,$$
where $\gamma>0$ is an absolute constant.
\endproclaim

From our proof more precise bounds
$$\split &\vl(\PP) \ \geq \ 
{\Gamma\left({m+n \over 2}\right)   \over 2 e^3 \sqrt{mn}  \pi^{m+n-2 \over 2} N^{m+n-1}}
{(mn)^{mn} \over (mn)!} \beta \\
&\qquad  \text{and} \\
&\vl(\PP) \ \leq \  {2e \lambda^{m+n-2} (mn)^{2m+2n-5/2} \over N^{m+n-1}}{(mn)^{mn} \over (mn)!} \beta,  \\ &\qquad \text{where} \\ &\lambda=\lambda(R,C)={n\over 2} \max_{i=1, \ldots, m} {N \over r_i}
 + {m \over  2}
\max_{j=1, \ldots, n} {N \over  c_j}\endsplit$$
follow.  When the margins are scaled,
 $(R,C) \longmapsto (tR, tC)$ for $t>0$, the volume of $\PP$ and both 
 the upper and the lower bounds get multiplied by  
 $t^{\dim \PP}$.

Computing $\beta$ reduces to finding the maximum of the concave 
function
$$f(X)=\sum_{ij} \ln x_{ij}$$ on the transportation polytope $\PP$ and hence can be
done efficiently (in time polynomial in $m+n$ and $\ln N$) by existing methods \cite{NN94}. 

Computing or estimating volumes of transportation polytopes has attracted considerable attention 
as a testing ground for methods of convex geometry \cite{Sch92}, combinatorics \cite{Pak00},
analysis and algebra \cite{BLV04}, \cite{BP03}, \cite{DLY03}. In a recent breakthrough \cite{CM07b}, Canfield 
and McKay obtained a precise asymptotic expression for the volume of the {\it Birkhoff polytope} 
(when $r_i=c_j=1$ for all $i$ and $j$) and in the more general case of all the row sums
being equal and all the column sums being equal. If $R=C=\left(1, \ldots, 1\right)$,
the formula of \cite{CM07b} gives 
$$\vl \PP={1 \over (2\pi)^{n-1/2} n^{(n-1)^2} }
\exp\left\{ {1 \over 3} + n^2 +O\left(n^{-1/2 +\epsilon} \right)\right\},$$
whereas the formula of Theorem 1.2 implies that, ignoring lower-order terms, we have
$$\vl \PP \approx {e^{n^2} \over n^{n^2}}$$
in that case (since by symmetry the maximum $\beta$ of the product of coordinates $x_{ij}$ 
 is attained at $x_{ij}=1/n$).
 Theorem 1.2 seems to be the only rigorously proven estimate of the volume 
of the transportation polytope available 
for general margins. 

We note that from the purely algorithmic perspective,
 volumes of polytopes and convex bodies can be computed in {\it randomized}
polynomial time, see \cite{Bol97} for a survey.

Theorem 1.1 can be extended to counting with weights.

Let us fix a non-negative matrix $W=\left(w_{ij}\right)$, which we call the matrix of {\it weights}.
We consider the following expression 
$$T(R, C; W)=\sum_{D=\left(d_{ij}\right)} \prod_{ij} w_{ij}^{d_{ij}},$$
where the sum is taken over all non-negative integer matrices $D$ with row sums $R$ 
and column sums $C$ and where we agree that $0^0=1$. For example, if 
$w_{ij} \in \{0,1\}$ for all $i,j$ then $T(R,C;W)$ is the number of $m \times n$ non-negative
integer matrices $D=\left(d_{ij}\right)$ with row sums $R$, column sums $C$ and such that
$d_{ij}=0$ whenever $w_{ij}=0$. This number can also be interpreted as the number 
of integer feasible flows in a bipartite graph
with vertices $u_1, \ldots, u_m$ and $v_1, \ldots, v_n$ and edges $(u_i, v_j)$ whenever 
$w_{ij}=1$ that satisfy the supply constraints $r_i$ at $u_i$ and the demand constraints $c_j$ at $v_j$.
Counting integer feasible flows in non-bipartite networks can be reduced to that for bipartite networks.
For example, if $w_{ij}=1$ for $j \leq i+1$ and $w_{ij}=0$ elsewhere, $T(R,C;W)$ is 
the Kostant partition function, see \cite{Ba07}, \cite{Ba08} for more examples and details.
We also note that $T(R,C; {\bold 1})=\#(R,C)$, where ${\bold 1}$ is the matrix of all 1's.

We prove the following extension of Theorem 1.1.

\proclaim{(1.3) Theorem} Let $R=\left(r_1, \ldots, r_m \right)$ and $C=\left(c_1, \ldots, c_n \right)$ 
be positive integer vectors such that $r_1 + \ldots + r_m =c_1 + \ldots + c_n =N$ and let 
$W=\left(w_{ij}\right)$ be an $m \times n$ non-negative matrix of weights.
Let us define a function
$$\split F(&\xx, \yy; W) =\left( \prod_{i=1}^m x_i^{-r_i} \right) \left(\prod_{j=1}^n y_j^{-c_j} \right)
\left( \prod_{ij} {1  \over 1- w_{ij} x_i y_j }\right) \\ & \text{for} \quad
\xx=\left(x_1, \ldots, x_m \right) \quad \text{and} \quad \yy=\left(y_1, \ldots, y_n \right)
\endsplit$$
and let
$$\rho=\rho(R,C; W)=\inf \Sb x_1, \ldots, x_m >0 \\ y_1, \ldots, y_n >0 \\ w_{ij} x_i y_j < 1
\ \text{for all} \  i,j \endSb F(\xx, \yy; W).$$

Then, for the number $T(R, C; W)$ of weighted non-negative integer matrices with 
row sums $r_1, \ldots, r_m$ and column sums $c_1, \ldots, c_n$, we have
$$\rho \ \geq \ T(R, C; W) \ \geq \ N^{-\gamma(m+n)} \rho,$$  
where $\gamma>0$ is an absolute constant.
\endproclaim
More precisely, the lower bound we prove is
$$\split T(R,C;W) \ \geq \  \ &{\Gamma\left({m+n \over 2} \right) \over 2e^5  
\pi^{{m+n-2 \over 2}} mn (N+mn)} \left({2\over (mn)^2 (N+1) (N+mn)}\right)^{m+n-1} \\
&\quad \times \left(\prod_{i=1}^m {r_i^{r_i} \over r_i!} \right) \left( \prod_{j=1}^n {c_j^{c_j} \over c_j!} \right) 
{N! (N+mn)! (mn)^{mn} \over  N^N  (N+mn)^{N+mn}(mn)!} \rho(R,C;W)\endsplit$$
provided $m+n \geq 10$.

As in Theorem 1.1, substituting $x_i=e^{-t_i}$ for $i=1, \ldots, m$ and 
$y_j=e^{-s_j}$ for $j=1, \ldots, n$ we reduce the problem of computing $\rho$ 
to the problem of finding the infimum of the convex function
$$\phi(\ttt,\sss)=\phi_{R,C}(\ttt,\sss)=\sum_{i=1}^m r_i t_i +\sum_{j=1}^n c_j s_j -
\sum_{ij} \ln\left(1-w_{ij} e^{-t_i -s_j} \right)$$
on the convex polyhedron
$$ s_i +t_j > \ln w_{ij} \quad \text{for all} \quad i,j.$$
Again, the value of $\rho$ can be computed efficiently, both in theory and in practice, by methods of convex optimization,
cf. \cite{NN94}.

For {\it positive} matrices $W=\left(w_{ij}\right)$ the infimum $\rho(R, C; W)$ in 
Theorem 1.3 is attained at a particular point and there is a convenient dual description 
of $\rho(R,C; W)$.

\proclaim{(1.4) Lemma}  Let $\PP=\PP(R,C)$ be the transportation polytope 
of the $m \times n$ non-negative matrices $X =\left(x_{ij}\right)$ with row sums $R$ and column sums $C$ and let  us fix an $m \times n$ positive matrix $W=\left(w_{ij}\right)$ of weights, so 
$w_{ij} >0$ for all $i,j$. For an $m \times n$ non-negative matrix $X=\left(x_{ij}\right)$ let us 
define
$$g(X; W) =\sum_{ij} \Bigl( \left(x_{ij} +1 \right) \ln \left(x_{ij}+1\right) -x_{ij} \ln x_{ij} + x_{ij} \ln w_{ij}
\Bigr).$$
Then $g(X; W)$ is a strictly concave function of $X$ and attains its maximum on 
$\PP$ at a unique positive matrix $Z=Z(R,C; W)$. One can write $Z=\left(z_{ij}\right)$ in the form
$$z_{ij}={w_{ij} \xi_i \eta_j  \over 1- w_{ij} \xi_i \eta_j}  \quad \text{for all} \quad i,j$$ 
and positive $\xi_1, \ldots, \xi_m; \eta_1, \ldots, \eta_n$ such that 
$w_{ij} \xi_i \eta_j < 1$ for all $i$ and $j$ and such that the infimum $\rho(R,C;W)$ in 
Theorem 1.3 is attained at $\xx^{\ast}=\left(\xi_1, \ldots, \xi_m \right)$ and 
$\yy^{\ast}=\left(\eta_1, \ldots, \eta_n \right)$:
$$\rho(R,C; W) =F\left(\xx^{\ast}, \yy^{\ast}; W \right).$$
Moreover, we have
$$\rho(R,C; W)=\exp\Bigl\{g(Z; W)\Bigr\}.$$
In particular, if $w_{ij}=1$ for all $i,j$, then 
$$\split &g(X)=g(X; {\bold 1}) = \sum_{ij} \Bigl( \left(x_{ij} +1 \right) \ln \left(x_{ij}+1\right) -x_{ij} \ln x_{ij}\Bigr) \quad \text{and}\\
&z_{ij}= {\xi_i \eta_j \over 1- \xi_i \eta_j} \quad \text{for all} \quad i,j, \endsplit$$
where $\xx^{\ast} =\left(\xi_1, \ldots, \xi_m \right)$ and $\yy^{\ast} =\left(\eta_1, \ldots, \eta_n \right)$
is a point where the minimum $\rho(R, C)=F\left(\xx^{\ast}, \yy^{\ast} \right)$ in Theorem 1.1 is attained.
Additionally,
$$\rho(R, C) =\exp\Bigl\{g(Z)\Bigr\}.$$
\endproclaim

The paper is structured as follows.

In Section 2, we consider consequences of Theorems 1.1 and 1.2 for the ``independence 
heuristic''. The heuristic was, apparently, first discussed by Good, see
\cite{Goo76}. It asserts that if we consider the space of non-negative 
integer $m \times n$ matrices with the total sum $N$ of entries as a probability space with the
 uniform measure then the event consisting of the matrices with the 
row sums $r_1, \ldots, r_m$ is ``almost independent''  from the event consisting 
of the matrices with the column sums $c_1, \ldots, c_n$. We show that if the row sums $r_i$ 
and the column sums $c_j$ are sufficiently generic then the independence heuristic 
tends to underestimate the number of tables as badly as within a factor of $\gamma^{mn}$ 
for some absolute constant $\gamma>1$. We see that in fact (rather counter-intuitively),
instead of independence, we have attraction (positive correlation) of the events.

In Section 3, we state a general result (Theorem 3.1), which provides a reasonably accurate 
estimate for the volume of the section of the standard simplex by a subspace of a  small 
codimension. Theorem 3.1 states that in a sufficiently generic situation the volume 
of the section is determined by the maximum value of the product of the coordinates 
of a point in the section. This estimate immediately implies Theorem 1.2 and is one of the two 
crucial ingredients in the proofs of Theorems 1.1 and 1.3. Theorem 3.1 appears to be 
new and may be interesting in its own right.

In Section 4, we state some preliminaries from convex geometry needed to prove Theorem 3.1.

In Section 5, we prove Theorems 3.1 and 1.2.

In Section 6, we describe the second main ingredient for the proofs of Theorems 1.1 and 1.3, 
the integral representation from \cite{Ba07} and \cite{Ba08} 
for the number $\#(R,C)$ of tables and the number $T(R, C; W)$ 
of weighted tables.  

In Section 7, we prove Theorems 1.1 and 1.3 and Lemma 1.4.
\bigskip
In what follows, we use $\gamma$ to denote a positive constant.

\head 2. The independence heuristic and the exponential attraction in the space of
matrices  \endhead 

\subhead (2.1) The independence heuristic \endsubhead  The following heuristic 
approach to counting contingency tables was suggested by Good \cite{Goo76}. Let us consider 
the space of all $m \times n$ non-negative integer matrices with the total sum of entries $N$ as 
a probability space with the uniform measure. Then the probability that a matrix from 
this space has row sums $R=\left(r_1, \ldots, r_m \right)$ is exactly 
$${N+mn-1 \choose mn-1}^{-1}\prod_{i=1}^m {r_i + n-1 \choose n-1}.$$
Similarly, the probability that a matrix has column sums $C=\left(c_1, \ldots, c_n \right)$
is exactly 
$${N+mn-1 \choose mn-1}^{-1}\prod_{j=1}^n {c_j + m-1 \choose m-1}.$$
Assuming that the two events are almost independent, one estimates the number 
$\#(R,C)$ of contingency tables by the {\it independence heuristic} $I(R,C)$:
$$I(R,C)={N+mn-1 \choose mn-1}^{-1}\prod_{i=1}^m {r_i + n-1 \choose n-1}\prod_{j=1}^n {c_j + m-1 \choose m-1}. \tag2.1.1$$
For example, if $m=n=4$, $R=(220, 215, 93, 64)$, $C=(108, 286, 71, 127)$ with $N=592$
then 
$$\#(R,C) =1225914276768514 \approx 1.226 \times 10^{15},$$
see \cite{DE85}, while
$$I(R,C) \approx 1.211 \times 10^{15}.$$

Given margins $R=\left(r_1, \ldots, r_n \right)$ and $C=\left(c_1, \ldots, c_m \right)$ such that
not all row sums $r_i$ are equal and not all column sums $c_j$ are equal,
we will construct a sequence of margins $(R_k, C_k)$, where $R_k$ is a $km$-vector and
$C_k$ is a $kn$-vector such that 
the ratio $\#(R_k,C_k)/I(R_k, C_k)$ grows as $\gamma^{k^2}$ for some $\gamma=\gamma(R,C)>1$.

\subhead (2.2) Cloning margins \endsubhead
Let us choose some margins $R=\left(r_1, \ldots, r_m \right)$ and $C=\left(c_1, \ldots, c_n\right)$
such that $r_1 + \ldots + r_m =c_1 + \ldots + c_n =N$.
For a positive integer $k$,  let us consider the new 
``clone'' margins 
$$\split &R_k=\left( \underbrace{kr_1, \ldots, kr_1}_{k \text{\ times}}, \ldots,  
\underbrace{kr_m, \ldots, kr_m}_{k \text{\ times}} \right)  \\
&C_k=\left( \underbrace{kc_1, \ldots, kc_1}_{k \text{\ times}}, \ldots,  
\underbrace{kc_n, \ldots, kc_n}_{k \text{\ times}} \right). 
\endsplit$$
In other words, we obtain margins $(R_k, C_k)$ if we choose an arbitrary matrix $X$ with row 
sums $R$ and column sums $C$, consider the $km \times kn$ block matrix $Y_k$ consisting 
of $k^2$ blocks $X$ and let $R_k$ be the row sums of $Y_k$ and let $C_k$ be the 
column sums of $Y_k$. Hence we consider $km \times kn$ matrices with the total sum of the matrix 
entries equal to $k^2 N$. 

One can check from the optimality condition (cf. Section 7.1) that if $\xx^{\ast}=\left(\xi_1, \ldots, \xi_m \right)$ and 
$\yy^{\ast}=\left(\eta_1, \ldots, \eta_n \right)$ is a point in Theorem 1.1 where the 
minimum $\rho(R, C)$ is attained then the minimum $\rho(R_k, C_k)$ is attained at 
the point 
$$\left( \underbrace{\xx^{\ast}, \ldots, \xx^{\ast}}_{\text{$k$ times}}, \underbrace{\yy^{\ast}, \ldots, 
\yy^{\ast}}_{\text{$k$ times}} \right).$$

Therefore, 
$$\rho(R_k, C_k) =\rho^{k^2}(R,C)$$ 
and by Theorem 1.1
$$\lim_{k \longrightarrow +\infty} \#(R_k, C_k)^{1/k^2} =\rho(R,C),$$
or, in other words,
$$\lim_{k \longrightarrow +\infty} {1 \over k^2} \ln \#(R_k, C_k) =\ln \rho(R,C). \tag2.2.1$$

Let us introduce the multivariate entropy function 
$$\HH\left(p_1, \ldots, p_d\right)=\sum_{i=1}^d p_i \ln {1 \over p_i},$$
where $p_1, \ldots, p_d$ are non-negative numbers such that $p_1 + \ldots + p_d=1$.
Using the standard asymptotic estimate for binomial coefficients (available, for example,
via Stirling's formula)
$$\lim_{k \longrightarrow + \infty} {1\over k} \ln {ka+ kb \choose ka} =(a+b) \ln (a+b) - a \ln a - b \ln b$$
 we deduce from (2.1.1) that 
$$\split  \lim_{k \longrightarrow +\infty}  {1 \over k^2} & \ln I(R_k, C_k)= \\
&-(N+mn) \HH\left({r_i + n \over N+mn},\ i=1, \ldots, m  \right) \\ 
&-(N+mn) \HH \left({c_j + m \over N + mn}, \ j=1, \ldots, n \right) \\
&-\sum_{i=1}^m r_i \ln r_i - \sum_{j=1}^n c_j \ln c_j  \\ &+ N \ln N + (N+mn) \ln (N +mn)
\endsplit \tag2.2.2$$

\subhead (2.3) The exponential attraction in the space of matrices \endsubhead
Let us choose margins $R=\left(r_1, \ldots, r_m \right)$ and $C=\left(c_1, \ldots, c_n \right)$ 
such that not all row sums $r_i$ are equal and not all column sums $c_j$ are equal.
Our goal is to show that 
$$\lim_{k \longrightarrow +\infty} {1 \over k^2} \ln \#(R_k, C_k)  \ > \ 
\lim_{k \longrightarrow +\infty} {1 \over k^2} \ln I(R_k, C_k), \tag2.3.1$$
so the ratio $\#(R_k, C_k)/I(R_k, C_k)$ grows as
$\gamma^{k^2}$ for some $\gamma=\gamma(R,C)>1$, as we clone margins 
$(R, C) \longmapsto (R_k, C_k)$.

By Lemma 1.4, we can write 
$$\ln \rho(R, C) \ = \ g(Z) \ \geq \ g(Y), \tag2.3.2$$
where $Y=\left(y_{ij}\right)$ is the {\it independence matrix} with $y_{ij}=r_i c_j/N$ for all $i,j$
and 
$$g(X) =\sum_{ij} \Bigl( \left(x_{ij}+1\right) \ln \left(x_{ij}+1\right) - x_{ij} \ln x_{ij} \Bigr).$$
On the other hand, it is easy to check that 
$$\split g(Y)
=&-(N+mn) \HH\left({r_i c_j + N \over N (N+mn)}, \quad i,j \right) \\
&\qquad \qquad -\sum_{i=1}^m r_i \ln r_i 
-\sum_{j=1}^n c_j \ln c_j \\ &\qquad \qquad +N\ln N + (N+mn) \ln (N+mn)  \endsplit \tag2.3.3$$
Let us consider the $m \times n$ matrix with the $(i,j)$-th entry equal to \newline
$(r_i c_j +N)/(N^2 +N mn)$. The $i$-th row sum of the matrix is $(r_i + n)/(N+mn)$, the $j$-th
column sum is $(c_j+m)/(N+mn)$ while the sum of all the entries of the matrix is 1.
Using the inequality relating the entropies of two partitions of a probability space with 
the entropy of the intersection of the partition (see, for example, \cite{Khi57}), we conclude
that 
$$\split &\HH\left({r_i c_j + N \over N (N+mn)}, \quad \matrix 1 \leq i \leq m \\ 1 \leq j \leq n \endmatrix\right)\\ &\qquad  \leq
  \HH\left({r_i + n \over N+mn},\ 1 \leq i \leq  m  \right)  \\  &\qquad \qquad +
 \HH \left({c_j + m \over N + mn}, \ 1 \leq j \leq n \right) \endsplit \tag2.3.4 $$
with the equality if and only if 
$${r_i + n \over N+mn} \cdot {c_j + m \over N + mn} ={r_i c_j + N \over N (N+mn)}
\quad \text{for all} \quad i,j. \tag2.3.5$$
Identities (2.3.5) are equivalent to $(N -r_i m)(N-c_j n)=0$, which, in turn, equivalent 
to all row sums being equal $r_i=N/m$ or all column sums being equal $c_j=N/n$.

Summarizing (2.2.1), (2.2.2), (2.3.2), and  (2.3.3) we conclude that inequality (2.3.1)
indeed holds if not all row sums $r_i$ are equal and not all column sums $c_j$ are equal.
 Therefore,
in the space of $km \times kn$ matrices with the sum $k^2 N$ of all entries the two events
$$\aligned & \RR_k: \quad \text{the vector of row sums of a matrix is} \quad R_k 
\\ &\qquad \qquad   \text{and}  \\
& \CC_k: \quad  \text{the vector of column sums of a matrix is} \quad C_k,
\endaligned \tag2.3.6
$$
instead of being asymptotically independent, attract exponentially in $k^2$, that is,
$${\Pr \left(\RR_k \cap  \CC_k \right) \over \left(\Pr \RR_k\right) \left(\Pr \CC_k \right)} 
\ \geq \ \gamma^{k^2}$$
for some $\gamma=\gamma(R,C)>1$ and all sufficiently large $k$.

Starting with non-constant margins $(R, C)$
 the cloning procedure $(R, C) \longmapsto (R_k, C_k)$ produces 
margins which stay away from from constant and maintain the density $N/mn$ separated 
from 0. Similar analysis shows that the phenomenon of attraction of the events $\RR_k$ and
$\CC_k$ defined by (2.3.6) holds for more general sequences of margins 
$(R_k, C_k)$ of growing dimensions which stay sufficiently away from uniform and sparse.

Two terms contribute to the difference $\ln \#(R, C) - \ln I(R,C)$:

first, the difference 
$g(Z)-g(Y)$, where $Z$ is the matrix of Lemma 1.4 at which the maximum 
of the function $g(X)=\sum_{ij} (x_{ij}+1) \ln (x_{ij}+1) - x_{ij} \ln x_{ij}$ on 
the transportation polytope $\PP(R, C)$ is attained  and $Y=\left(y_{ij}\right)$ is the 
independence matrix $y_{ij}=r_i c_j/N$, cf. (2.3.2);

 and second, the difference (multiplied by $(N+mn)$) between the entropies 
on the right hand side of (2.3.4) and the left hand side of (2.3.4). 

As long as either of these 
differences remains large enough to overcome the error term of 
$O\bigl((m+n) \ln N\bigr)$ coming from Theorem 1.1, we 
have the asymptotic positive correlation of sequences of events $\RR_k$ and $\CC_k$ in (2.3.6).

On the other hand, the independence estimate $I(R,C)$ produces a reasonable approximation
to $\#(R,C)$ in the cases of sparse tables  (cf. \cite{GM07}) and tables with constant 
margins (cf. \cite{CM07a}). One can show that if all row sums are equal or if all column 
sums are equal then indeed
$$\lim_{k \longrightarrow +\infty} {1 \over k^2} \ln \#(R_k, C_k)  \ = \ 
\lim_{k \longrightarrow +\infty} {1 \over k^2} \ln I(R_k, C_k), $$
where $(R_k, C_k)$ are cloned margins $(R,C)$.
Indeed, if all $r_i$ are equal then the symmetry argument shows that the matrix 
$Z=\left(z_{ij}\right)$ in Lemma 1.4 satisfies $z_{ij}=c_j/m$ for all $i$ and $j$, and, similarly, if all $c_j$ 
are equal then we have $z_{ij}=r_i/n$ for all $i,j$. In either case we have $Z=Y$ in (2.3.2) and, as 
we have already discussed, equations (2.3.5) hold as well.

\head 3. The volume of a section of a simplex \endhead

Let $\AA$ be the affine hyperplane in ${\Bbb R}^d$ defined by the equation
$$\sum_{i=1}^d x_i =1$$
and let $\Delta \subset \AA$ be the standard $(d-1)$-dimensional open simplex defined by 
the inequalities 
$$x_i > 0 \quad \text{for} \quad i=1, \ldots, d.$$
We consider the Euclidean structure in $\AA$ induced from ${\Bbb R}^d$.
In particular, if $K \subset \AA$ is an $m$-dimensional convex body, by $\vl_m(K)$
we denote the $m$-dimensional volume of $K$ with respect to that Euclidean structure.
For $m=d-1$ we denote $\vl_m$ just by $\vl$. In particular,
$$\vl (\Delta)={\sqrt{d} \over (d-1)!}.$$

Let $L \subset \AA$ be an affine subspace intersecting  $\Delta$. 
Suppose that $\dim L=d-k-1$, so the the codimension of $L$ in $\AA$ is $k \geq 1$.
 Our aim is to estimate the volume of the intersection
$\vl_{d-k-1} (L \cap \Delta)$ within a reasonable accuracy when the codimension $k$ of 
$L$ is small. It turns out that the volume is controlled by one particular quantity, namely 
the maximum value of the product of the coordinates of a point $x \in \Delta \cap L$.

Our result is as follows.
\proclaim{(3.1) Theorem} Let $\AA \subset {\Bbb R}^d$ be the affine hyperplane 
defined by the equation $x_1 + \ldots + x_d=1$ and let $\Delta \subset \AA$ be the 
standard $(d-1)$-dimensional open simplex defined by the inequalities $x_1 > 0, \ldots, x_d >  0$.
 
Let $L \subset \AA$ be an affine subspace intersecting 
$\Delta$ and such that $\dim L=d-k-1$ where $k \geq 1$. Suppose that the maximum of the function
$$f(x)=\sum_{i=1}^d \ln x_i$$ 
on $\Delta \cap L$ is attained at $a=\left(\alpha_1, \ldots, \alpha_d \right)$.
\roster
\item We have 
$${\vl_{d-k-1}(\Delta \cap L) \over \vl(\Delta)} \ \geq \ \gamma {1 \over d^2 \omega_k} 
d^d e^{f(a)},$$
where 
$$\omega_k ={\pi^{k/2} \over \Gamma(k/2+1)}$$ is the volume of the $k$-dimensional 
unit ball and $\gamma>0$ is an absolute constant (one can choose $\gamma=1/2e^3 \approx 
0.025$).
\item Suppose that 
$$ \alpha_i \geq {\epsilon \over d} \quad \text{for some} \quad \epsilon>0 
\quad \text{and} \quad i=1, \ldots, d.$$
Then 
$${\vl_{d-k-1}\left(\Delta \cap L \right) \over \vl(\Delta)} \ \leq\  \gamma\left({ d^2 \over 2 \epsilon}\right)^k 
d^d e^{f(a)},$$
where $\gamma >0$ is an absolute constant (one can choose $\gamma=2e \approx 5.44$).
\endroster
\endproclaim

We are interested in the situation of $k \sim \sqrt{d}$, so ignoring lower-order terms in 
the logarithmic
order, we get
$$\vl_{d-k-1} (\Delta \cap L) \sim {d^d \over d!} e^{f(a)} \sim e^d \prod_{i=1}^d \alpha_i,$$
provided the maximum value of the product of the coordinates of a point $x \in \Delta \cap L$
is attained at $a=\left(\alpha_1, \ldots, \alpha_d \right)$ and all $\alpha_i$ are not too small.

Let 
$$c=\left({1 \over d}, \ldots, {1 \over d} \right)$$
be the center of the simplex $\Delta$.

We deduce Theorem 3.1 from the following result.
\proclaim{(3.2) Theorem} Let $\AA \subset {\Bbb R}^d$ be the affine hyperplane 
defined by the equation $x_1 + \ldots + x_d=1$ and let $\Delta \subset \AA$ be the 
standard $(d-1)$-dimensional open simplex defined by the inequalities $x_1 > 0, \ldots, x_d >0$.

Let $H \subset \AA$ be an affine hyperplane in $\AA$ intersecting $\Delta$. If $H$ does not pass through the center $c$ of $\Delta$, let $H^- \subset \AA$ be 
the open halfspace bounded by $H$ that does not contain $c$ and if $H$ passes 
through $c$ let $H^- \subset \AA$ be either of the open halfspaces bounded by $H$.
 
 Suppose that the function
 $$f(x)=\sum_{i=1}^d \ln x_i \quad \text{where} \quad x=\left(x_1, \ldots, x_d\right)$$
 attains its maximum on $\Delta \cap H$ at a point $a=\left(\alpha_1, \ldots, \alpha_d\right)$.

Then,
for some absolute constant $\gamma>0$ we have 
$$d^d e^{f(a)}  \geq \ {\vl(\Delta \cap H^-) \over \vl(\Delta)} \geq {\gamma \over d^2}  d^d e^{f(a)}.$$ 
We can choose $\gamma=1/2e^3 \approx 0.025$.
\endproclaim

\head 4. Preliminaries from convex geometry \endhead

We recall that $\AA \subset {\Bbb R}^d$ is the affine
hyperplane defined by the equation $x_1 + \ldots + x_d=1$, that 
$\Delta \subset \AA$ is the open simplex defined by the inequalities $x_i > 0$ for 
$i=1, \ldots, d$, and that $c=\left(1/d, \ldots, 1/d \right)$ is the center of $\Delta$.
We need some results 
regarding central hyperplane sections of $\Delta$.

\proclaim{(4.1) Lemma} 
Let $H \subset \AA$ be an affine hyperplane in $\AA$ passing through the center $c$ 
of $\Delta$. 
\roster
\item Let $H^+$ and $H^{-}$ be the open halfspaces bounded 
by $H$. Then for some absolute constant $\gamma>0$ we have
$${\vl \left(\Delta \cap H^+ \right) \over \vl(\Delta)} \geq \gamma \quad  \text{and} \quad
{\vl \left(\Delta \cap H^- \right) \over \vl(\Delta)} \geq \gamma.$$
One can choose $\gamma=1/e \approx 0.37$.
\item For some absolute constant $\gamma>0$ we have
$$ {\vl_{d-2}(\Delta \cap H) \over \vl (\Delta)} \geq \gamma.$$
One can choose $\gamma=1/2e \approx 0.18$.
\endroster
\endproclaim
\demo{Proof} Part (1) is a particular case of a more general result of Gr\"unbaum \cite{Gr\"u60}
on hyperplane sections through the centroid of a convex body.
 In fact, in dimension $d$ one can choose 
 $$\gamma_d=\left(1 -{1 \over d} \right)^{d-1} > {1 \over e}.$$

As K. Ball and M. Fradelizi explained to the author, a stronger estimate than that of Part (2) 
can be obtained by combining techniques of \cite{Bal88} and \cite{Frad97}. Nevertheless,
we present a proof of Part (2) below since the same approach is used later 
in the proof of Theorem 3.1.

To prove Part (2), let $H^{\bot} \subset \AA$ be a line orthogonal to $H$. Let us consider
the orthogonal projection $pr: \AA \longrightarrow H^{\bot}$ and let $Q=pr(\Delta)$ be
the image of the simplex. Since $\Delta$ is contained in a ball of radius 1, $Q$ is an interval 
of length at most 2. 

Let $y_0=pr(H)$ and for $y \in Q$ let 
$$\nu(y)=\vl_{d-2} \left(pr^{-1}(y)\right)$$
be the volume of the inverse image of $y$. By the Brunn-Minkowski inequality, the 
function $\nu$ is log-concave, see \cite{Bal88}, \cite{Bal97}.

Our goal is to bound $\nu(y_0)$ from below. The point $y_0$ splits the interval
$Q$ into two subintervals, $Q^+=pr(\Delta \cap H^+)$ and 
$Q^-=pr(\Delta \cap H^-)$ of length at most 2 each.

We have
$$\int_{Q^+} \nu(y) \ dy =\vl\left(\Delta \cap H^+\right) \quad \text{and} \quad 
\int_{Q^-} \nu(y) \ dy = \vl\left(\Delta \cap H^-\right).$$
Using Part (1) we conclude that there exist $y^+ \in Q^+$ and $y^- \in Q^-$ such that
$${\nu(y^+) \over \vl (\Delta)} \geq {1 \over 2e} \quad \text{and} \quad
{\nu(y^-) \over \vl (\Delta)} \geq {1 \over 2e}.$$

Since $y_0$ is a convex combination of $y^+$ and $y^-$, by  the log-concavity of $\nu$ we must have
$${\nu(y_0) \over \vl (\Delta)} \geq {1 \over 2e},$$
as desired.
{\hfill \hfill \hfill} \qed
\enddemo

Let us choose a point $a=\left(\alpha_1, \ldots, \alpha_d\right)$ in $\Delta$,  and 
let us consider the {\it projective transformation} $T_a: \Delta \longrightarrow \Delta$
$$\split T_a(&x)=y \quad \text{where} \quad y_i={\alpha_i x_i \over \alpha_1 x_1 + \ldots +\alpha_d x_d}
\quad \text{for} \\ &x=\left(x_1, \ldots, x_d \right) \quad \text{and} \quad y=\left(y_1, \ldots, y_d\right). \endsplit$$
The inverse transformation is $T_b$ for $b=\left(\alpha_1^{-1}, \ldots, \alpha_d^{-1} \right)$.
Clearly,
$$T_a(c)=a,$$
where $c$ is the center of $\Delta$.
For $x \in \Delta$, the derivative $DT_a(x)$ is a linear transformation
$$DT_a(x): \quad \HHH \longrightarrow \HHH,$$
where $\HHH$ is the hyperplane $x_1+ \ldots + x_d=0$ in ${\Bbb R}^d$.

Our immediate goal is to compute the Jacobian $|DT_a(x)|$ at $x \in \Delta$.

\proclaim{(4.2) Lemma} Let us choose a point $a=\left(\alpha_1, \ldots, \alpha_d \right)$ 
in the simplex $\Delta$ and let us consider the projective transformation
$$T_a: \Delta \longrightarrow \Delta$$
defined by the formula
$$T_a(x)=y \quad \text{where} \quad y_i={\alpha_i x_i \over \alpha_1 x_1 + \ldots + \alpha_d x_d}$$
for $x=(x_1, \ldots, x_d)$ and $y=(y_1, \ldots, y_d)$.

Let $DT_a(x): \HHH \longrightarrow \HHH$ be the derivative of $T_a$ at $x \in \Delta$ and
$|DT_a(x)|$ the corresponding value of the Jacobian. 
Then 
$$|DT_a(x)|= {\alpha_1 \cdots \alpha_d \over (\alpha_1 x_1 + \ldots + \alpha_d x_d)^d}.$$
\endproclaim
\demo{Proof} Let us consider $T_a$ as defined in a neighborhood of $x$ in ${\Bbb R}^d$ with 
values in ${\Bbb R}^d$ and 
let $DT_a(x)$ be the $d \times d$ matrix of the derivative 
$$DT_a(x)=\left({\partial y_i \over \partial x_j} \right)$$
in the standard basis of ${\Bbb R}^d$. Then the $i$-th diagonal entry of $DT_a(x)$
is 
$${\alpha_i \over (\alpha_1 x_1 + \ldots + \alpha_d x_d)} - 
{\alpha_i^2 x_i \over (\alpha_1 x_1 + \ldots + \alpha_d x_d)^2},$$
while the $(i,j)$-th entry for $i \ne j$ is 
$${-\alpha_i \alpha_j x_i \over (\alpha_1 x_1 + \ldots + \alpha_d x_d)^2}.$$ 

Let 
$$\beta={1 \over \alpha_1 x_1 + \ldots + \alpha_d x_d},$$ let $B$ be the diagonal matrix
with the diagonal entries $\alpha_1, \ldots, \alpha_d$, and let $C$ be the matrix with the 
$(i,j)$th entry equal to $\alpha_i \alpha_j x_i$ for all $1 \leq i,j \leq d$. Then we can write 
$$DT_a(x)=\beta B-\beta^2 C=\beta(B-\beta C).$$

Since $DT_a(x)$ maps ${\Bbb R}^d$ onto $\HHH$ and $\HHH$ is an invariant subspace of 
$DT_a(x)$, the value of the Jacobian we are interested in is the product of the $(d-1)$ non-zero 
eigenvalues of $DT_a(x)$ (counting algebraic multiplicities), which is equal to the $(d-1)$-st elementary symmetric function in 
the eigenvalues of $DT_a(x)$, which is equal to the sum of the $d$ of the $(d-1) \times (d-1)$ 
principle minors of $DT_a(x)$.

Let $B_i$ and $C_i$ be the $(d-1) \times (d-1)$ matrices obtained from $B$ and $C$ respectively 
by crossing out the $i$th row and column.

Hence
$$B_i -\beta C_i = B_i (I - \beta B_i^{-1} C_i),$$
where $I$ is the $(d-1) \times (d-1)$ identity matrix. Now $B_i^{-1} C_i$ is a matrix of rank 1 
with the non-zero eigenvalue equal to the trace of $B_i^{-1} C_i$, which is 
$$\sum_{j \ne i} \alpha_j x_j.$$
Hence 
$$\det\left(I - \beta B_i^{-1} C_i \right)= 1- \beta \sum_{j \ne i} \alpha_j x_j =
{\alpha_i x_i \over \alpha_1 x_1 + \ldots + \alpha_d x_d}$$
and 
$$\det B_i (I- \beta B_i^{-1} C_i ) = {\alpha_1 \cdots \alpha_d x_i \over \alpha_1 x_1 + \ldots 
+ \alpha_d x_d}.$$
Therefore, the sum of the $d$ of $(d-1) \times (d-1)$ principle minors of $B-\beta C$ is 
$${\alpha_1 \ldots \alpha_d \over \alpha_1 x_1 + \ldots + \alpha_d x_d}$$
and the sum of the $(d-1)\times (d-1)$ principle minors of $DT_a(x)=\beta(B-\beta C)$ is 
$${\alpha_1 \ldots \alpha_d \over (\alpha_1 x_1 + \ldots + \alpha_d x_d)^d},$$
as desired.
{\hfill \hfill \hfill} \qed
\enddemo

Next, we will need a technical estimate, which shows that if the volume of the section of 
the simplex by an affine subspace of a small codimension is sufficiently large and if the 
subspace cuts sufficiently deep into the simplex then a neighborhood of the section in the 
simplex has a sufficiently large volume.

\proclaim{(4.3) Lemma} Let $L \subset \AA$ be an affine subspace,
$\dim L=d-k-1$. Suppose that there is a point $a \in L \cap \Delta$,
$a=\left(\alpha_1, \ldots, \alpha_d\right)$ such that
$$\alpha_i \geq {\epsilon \over d} \quad \text{for} \quad i=1, \ldots, d$$
and some $\epsilon >0$.

Let
$$\|x\|_{\infty}=\max\bigl\{|x_i|  \quad \text{for} \quad i=1, \ldots, d \bigr\} \quad 
\text{for} \quad i=1, \ldots, d$$ 
 and let us define a neighborhood $Q$ of 
$\Delta \cap L$ by 
$$Q=\Bigl\{x \in \Delta: \quad \|x -y\|_{\infty} \leq {\epsilon \over d^2} \quad 
\text{for some} \quad y \in \Delta \cap L \Bigr\}.$$
Then, for any affine hyperplane $H \subset \AA$ passing through $L$ we have
$$\vl ( Q \cap H^+), \quad \vl(Q \cap H^-) \geq \gamma \left({2 \epsilon \over d^2} \right)^k
\vl_{d-k-1} (\Delta \cap L),$$
where $H^+$ and $H^-$ are the halfspaces bounded by $H$ and $\gamma>0$ is an 
absolute constant. One can choose $\gamma=1/2e \approx 0.18$.
\endproclaim
\demo{Proof} Let
$$Q_0 =\left\{{1 \over d} a + {d-1 \over d} x:\quad x \in \Delta \cap L \right\}.$$
Since $Q_0$ is the contraction of $\Delta \cap L$ we have
$$\vl_{d-k-1} Q_0 =\left({d-1 \over d}\right)^{d-k-1} \vl_{d-k-1} (\Delta \cap L) \geq {1 \over e} \vl_{d-k-1} (\Delta \cap L).$$

Moreover, for any $x \in Q_0$, 
$x=\left(x_1, \ldots, x_d \right)$, we have 
$$x_i \geq {\epsilon \over d^2} \quad \text{for} \quad i=1, \ldots, d.$$
For every point $x \in Q_0$ let us consider the cube
$$I_x =\Bigl\{y \in {\Bbb R}^d:  \quad \|y-x\|_{\infty} \leq { \epsilon \over d^2} \Bigr\}.$$
Then $\left(I_x \cap \AA\right) \subset \Delta$. The intersection of $I_x$ with 
the $k$-dimensional affine subspace 
$L^{\bot}_x \subset \AA$ orthogonal to $L$ and passing through $x$ is centrally symmetric 
with respect to $x$ and, by Vaaler's Theorem \cite{Vaa79}, satisfies
$$\vl_k \left(I_x \cap L^{\bot}_x \right) \geq \left({2 \epsilon \over d^2} \right)^k.$$
The proof now follows.
{\hfill \hfill \hfill} \qed
\enddemo

\head 5. Proofs of Theorems 1.2, 3.1, and 3.2\endhead

We prove Theorem 3.2 first.

\subhead (5.1) Proof of Theorem 3.2 \endsubhead 
If $c \in H$ the result follows by Lemma 4.1. Hence we assume that $c \notin H$.

The hyperplane $H$ is orthogonal to the gradient of $f(x)$ at $x=a$ and passes through $a$, from which it follows
that $H$ can be defined in $\AA$ by the equation
$$\sum_{i=1}^d {x_i \over \alpha_i}  =d,$$
while the halfspace $H^-$ is defined by the inequality 
$$\sum_{i=1}^d {x_i \over \alpha_i} <  d.$$

Let us consider the projective transformation $T_a: \Delta \longrightarrow \Delta$ 
defined by the formula of Lemma 4.2. Hence $T_a(c)=a$. Moreover, 
the inverse image $T^{-1}_a(H)$ is the hyperplane $H_0$ defined in $\AA$ by the equation 
$$\sum_{i=1}^d \alpha_i x_i={1 \over d}$$
and the inverse image $T^{-1}_a\left(\Delta \cap H^-\right)$ is the 
intersection $\Delta \cap H_0^-$, where
 $H_0^-$ is the halfspace defined by the inequality 
$$\sum_{i=1}^d \alpha_i x_i > {1 \over d}.$$
By Lemma 4.2, we have 
$$|DT_a(x)|= {\alpha_1 \cdots \alpha_d \over \left( \alpha_1 x_1 + \ldots + \alpha_d x_d \right)^d} 
< d^d \prod_{i=1}^d \alpha_i \quad \text{for all} \quad x \in \Delta \cap H_0^-.$$
Since
$$\vl(\Delta \cap H^-) =\int_{\Delta \cap H_0^-} |DT_a(x)| \  dx, \tag5.1.1$$
the upper bound follows. 

Let us prove the lower bound. By Part (2) of Lemma 4.1, 
$${\vl_{d-2} \left(\Delta \cap H_0\right) \over \vl( \Delta)} \geq {1 \over 2e}.$$
We recall that $H_0$ passes through the center 
of the simplex and apply Lemma 4.3 with $\epsilon=1$.
Namely, we define 
$$Q=\Bigl\{x \in \Delta: \quad \|x-y\|_{\infty} \leq {1 \over d^2} \quad \text{for some} \quad 
y \in \Delta \cap H_0 \Bigr\}$$
and conclude that by Lemma 4.3
$$\vl (Q \cap H_0^-) \  \geq \ \left({1 \over 2e} \right) \left({2 \over d^2}\right) \vl_{d-2} \left(\Delta \cap H_0 \right) 
\ \geq \ {1 \over 2 e^2 d^2}  \vl(\Delta).$$
We note that for every $x \in Q$ we have 
$$\sum_{i=1}^d \alpha_i x_i \leq {1 \over d} +{1 \over d^2} ={d+1 \over d^2}.$$
By (5.1.1)
$$\split \vl (\Delta \cap H^-) \geq &\int_{Q \cap H_0^-} |DT_a(x)| \ dx \geq 
\left({d^2 \over d+1}\right)^d \vl(Q \cap H_0^-) \prod_{i=1}^d \alpha_i \\ \geq &{1 \over 2 e^3} {1 \over d^2}
 d^d \vl(\Delta) \prod_{i=1}^d \alpha_i  ,\endsplit $$
which completes the proof.
{\hfill \hfill \hfill} \qed

Next, we prove Theorem 3.1.

\subhead (5.2) Proof of Theorem 3.1 \endsubhead
The proof of Part (1) is similar to that of Part (2) of Lemma 4.1. Let $L^{\bot} \subset \AA$ 
be a $k$-dimensional subspace orthogonal to $L$ in $\AA$ and let 
$$pr:  \AA \longrightarrow L^{\bot}$$ be the orthogonal projection. Let $Q \subset L^{\bot}$,
$Q=pr(\Delta)$, 
be the image of the simplex. Clearly, $Q$ lies in a ball of radius $1$, so 
$$\vl_k Q \ \leq \ \omega_k.$$ 

For $y \in Q$, let 
$$\nu(y)=\vl_{d-k-1}\left(pr^{-1}(y)\right)$$
be the volume of the inverse image of $y$. By the Brunn-Minkowski inequality, the 
function $\nu$ is log-concave, so for every $\alpha >0$ the set 
$$\Bigl\{ y \in Q: \quad \nu(y) \geq \alpha \Bigr\}$$ is convex.
Moreover, for all Borel sets $Y \subset Q$ we have 
$$\int_Y \nu(y) \ dy = \vl\left(pr^{-1}(Y) \right).$$
We want to estimate $\nu(y_0)$ for $y_0=pr(L)$. 
Let $H \subset L^{\bot}$ be an affine hyperplane through $y_0$ and let $H^+, H^- \subset L^{\bot}$ 
be open halfspaces bounded by $H$. Then 
$\tilde{H}=pr^{-1}(H)$ is an affine hyperplane in $\AA$ containing $L$ and 
 $pr^{-1}(H^-)$ and $pr^{-1}(H^+)$ are the corresponding open halfspaces of $\AA$ 
bounded by $\tilde{H}$.

Since the maximum value of $f$ on $\Delta \cap \tilde{H}$ is at least as big as 
the maximum value of $f$ on $\Delta \cap L$,  by Theorem 3.2 we have 
$$\vl pr^{-1}\left(H^{\pm} \cap Q \right) \geq  {1 \over 2 e^3} {1 \over d^2} 
d^d e^{f(a)} \vl(\Delta).$$
Since 
$$\vl pr^{-1} \left(H^{\pm} \cap Q  \right)=\int_{H^{\pm} \cap Q} \nu(y) \ dy,$$
We conclude that there exist points $y^+ \in H^+$ and $y^- \in H^{-}$ such that
$$\nu(y^+), \ \nu(y^-) \ \geq\  {1 \over 2e^3} {1 \over d^2} d^d e^{f(a)} {\vl (\Delta) \over \vl_k Q}
\ \geq \  {1 \over 2e^3} {1 \over d^2} d^d e^{f(a)} {\vl (\Delta) \over \omega_k}. \tag5.2.1$$ 
In other words, for any affine hyperplane $H \subset L^{\bot}$ through $y_0$ on either 
side of the hyperplane there are points $y^+, y^-$ for which inequality (5.2.1) holds.
Hence $y_0$ lies in the convex hull of points $y$ for which the inequality holds.
The proof of Part (1) follows by the log-concavity of $\nu$.

Let us prove Part (2).
 Since $a$ is the maximum point of the strictly concave function
$$f(x)=\sum_{i=1}^d \ln x_i $$
on $\Delta \cap L$, the gradient of $f$ at $a$ is orthogonal to $L$. Hence $L$ is orthogonal 
to the vector
$$\left({1 \over \alpha_1}, \ldots, {1 \over \alpha_d} \right).$$
If $a \ne c$, let $H \subset \AA$ be the affine hyperplane defined by the equation
$$\sum_{i=1}^d {x_i \over \alpha_i} =d$$
and if $a=c$ let $H$ be any affine hyperplane containing $L$.
In either case $L \subset H$ and the maximum values of $f$ on $\Delta \cap H$ and on
 $\Delta \cap L$ coincide and are equal to $f(a)$.
Therefore, by Theorem 3.2, we have
$$\vl(\Delta \cap H^-) \ \leq \ d^d e^{f(a)} \vl(\Delta) \tag5.2.2$$
for some open halfspace $H^-$ bounded by $H$.

We apply Lemma 4.3. Namely, we let 
$$Q=\Bigl\{x \in \Delta: \quad \|x-y\|_{\infty} \leq {\epsilon \over d^2} \quad \text{for some}
\quad y \in \Delta \cap L \Bigr\}.$$
Then, by Lemma 4.3,
$$\vl (Q \cap H^-) \ \geq \ {1 \over 2e}  \left({2 \epsilon \over d^2} \right)^k \vl_{d-k-1}(\Delta \cap L).$$
Since 
$$\vl(Q \cap H^-)\ \leq \ \vl(\Delta \cap H^-),$$ 
we get the upper bound from (5.2.2).
{\hfill \hfill \hfill} \qed

Finally, we prove Theorem 1.2.
\subhead (5.3) Proof of Theorem 1.2 \endsubhead Let us consider the contracted polytope
$N^{-1} \PP$ defined by the equations
$$\sum_{j=1}^n x_{ij} ={r_i \over N} \quad \text{for} \quad i=1, \ldots, m, \qquad
\sum_{i=1}^m x_{ij} ={c_j  \over N}\quad \text{for} \quad j=1, \ldots, n$$
and inequalities
$$x_{ij} > 0 \quad \text{for all} \quad i,j.$$
Then $N^{-1}\PP$ can be represented as an intersection of the standard simplex in the 
space of $m \times n$ matrices and an affine subspace of dimension $(m-1)(n-1)$. 
We are going to use Theorem 3.1. Let $A=\left(\alpha_{ij}\right)$, $A \in N^{-1}\PP$, be the 
point maximizing the product of the coordinates. Writing the optimality condition 
for 
$$f(X)=\sum_{ij} \ln x_{ij}$$
on $N^{-1}\PP$, 
we conclude that 
$${1 \over \alpha_{ij}} =\lambda_i +\mu_j \quad \text{for all} \quad i,j$$
and some $\lambda_1, \ldots \lambda_m$ and $\mu_1, \ldots, \mu_n$.
Since $\lambda_i +\mu_j >0$ for all $i,j$, we may assume that $\lambda_i, \mu_j>0$ for all $i,j$. If $\lambda_i > nN/r_i$ for some 
$i$ then $ \alpha_{ij}< r_i/nN$ for all $j$, which is a contradiction. If $\mu_j > mN/c_j$ for 
some $j$ then $\alpha_{ij}< c_j/mN$ for all $i$ which is a contradiction. Hence
$\lambda_i \leq nN/r_i$ for $i=1, \ldots, m$ and $\mu_j \leq mN/c_j$ for $j=1, \ldots, n$,
from which 
$$\alpha_{ij} \geq {r_i c_j \over Nnc_j +Nmr_i} = {1 \over (nN/r_i) +(mN/c_j)}
 \quad \text{for all} \quad i,j.$$
The proof now follows by Theorem 3.1 with $d=mn$, $k=m+n-2$, and 
$$\epsilon =\left(n \max_{i=1, \ldots, m} {N \over r_i} + m \max_{j=1, \ldots, n} {N \over c_j} \right)^{-1}.$$
{\hfill \hfill \hfill} \qed

\head 6. An integral representation for the number of contingency tables \endhead

In this section, we recall bounds for $\#(R,C)$ obtained in \cite{Ba07} and \cite{Ba08}.
\subhead (6.1) Matrix scaling \endsubhead Our estimates for the number $\#(R,C)$ of contingency 
tables essentially use the theory of {\it matrix scaling}, see \cite{Si64}, \cite{MO68}, \cite{RS89}.
Let us fix non-negative vectors $R=\left(r_1, \ldots, r_m \right)$, $C=\left(c_1, \ldots, c_n \right)$, 
such that 
$$\sum_{i=1}^m r_i =\sum_{j=1}^n c_j=N.$$
Then for every $m \times n$ positive matrix $X=\left(x_{ij}\right)$ there exist
a positive $m \times n$ matrix $L=\left(l_{ij}\right)$ and positive numbers 
$\lambda_1, \ldots, \lambda_m$ and 
$\mu_1, \ldots, \mu_n$
such that 
$$\split &\sum_{j=1}^n l_{ij}=r_i  \quad \text{for} \quad i=1, \ldots, m, \\
&\sum_{i=1}^m l_{ij}=c_j \quad \text{for} \quad j=1, \ldots, n, \qquad \text{and} \\
&x_{ij}=\lambda_i \mu_j l_{ij} \quad \text{for all} \quad i,j. \endsplit  \tag6.1.1$$
Moreover, given $X$, the matrix $L$ is unique while the numbers 
$\lambda_i$ and $\mu_j$ are unique up to a re-scaling:
$$\split &\lambda_i \longmapsto \lambda_i \tau \quad \text{for} \quad i=1, \ldots, m \\
&\mu_j \longmapsto \mu_j \tau^{-1} \quad \text{for} \quad j=1, \ldots, n \endsplit$$
and some $\tau>0$.

\subhead (6.2) Function $\phi$ \endsubhead
This allows us to define a function
$$\phi(X)=\phi_{R,C}(X)=\left(\prod_{i=1}^m \lambda_i^{r_i} \right)
\left( \prod_{j=1}^n \mu_j^{c_j} \right),$$
where $\lambda_i$ and $\mu_j$ are numbers such that equations (6.1.1) hold,
on positive $m \times n$ matrices $X$.
It turns out that $\phi$ is continuous (it is also log-concave but we don't use that), positive homogeneous
of degree $N$, $$\phi(\alpha X)=\alpha^N \phi(X)$$
for $\alpha>0$ and positive matrix $X$, and monotone
$$\phi(X) \geq \phi(Y)$$
provided $X$ and $Y$ are positive matrices satisfying $x_{ij} \geq y_{ij}$ for all $i,j$, see,
for example, \cite{Ba07} and \cite{Ba08}.

Alternatively, $\phi(X)$ can be defined by 
$$\phi(X)=\min_{a,b} \left({1 \over N} \sum_{ij} x_{ij} \alpha_i \beta_j \right)^N,$$
where the minimum is taken over all positive $m$-vectors $a=\left(\alpha_1, \ldots, \alpha_m\right)$
and positive $n$-vectors $b=\left(\beta_1, \ldots, \beta_n \right)$ 
satisfying 
$$\prod_{i=1}^m \alpha_i^{r_i}=\prod_{j=1}^n \beta_j^{c_j}=1,$$
see also \cite{MO68}.

\subhead (6.3) The bounds \endsubhead
Let us identify the space of $m \times n$ matrices with  Euclidean space ${\Bbb R}^d$ 
for $d=mn$, let $\AA \subset {\Bbb R}^d$ be the affine hyperplane defined by the 
equation 
$$\sum_{ij} x_{ij} =1,$$
and let $\Delta \subset \AA$ be the standard open simplex defined by the inequalities 
$$x_{ij} > 0 \quad \text{for all} \quad i,j$$
with the Lebesgue measure $dX$ induced from the Euclidean structure in 
${\Bbb R}^d$.
It is proved in \cite{Ba07} and \cite{Ba08} that 
$$\aligned \#(R,C) \ \geq \ {N! (N+mn-1)! \over N^N \sqrt{mn}}  &\left(\prod_{i=1}^m {r_i^{r_i} \over r_i!}\right)
\left( \prod_{j=1}^n {c_j^{c_j} \over c_j!} \right) \int_{\Delta} \phi(X) \ dX  \\
&\qquad \qquad \text{and} \\
\#(R,C) \ \leq\  {(N+mn-1)! \over \sqrt{mn}} &\min\left\{ \prod_{i=1}^m {r_i^{r_i} \over r_i!}, \quad
\prod_{j=1}^n {c_j^{c_j} \over c_j!}\right\} \int_{\Delta} \phi(X) \ dX.  \endaligned \tag6.3.1 $$
Therefore, we have an approximation within up to an $N^{\gamma(m+n)}$ factor
for some absolute constant $\gamma>0$:
$$\#(R, C) \approx e^N  (N+mn)!  \int_{\Delta} \phi(X) \ dX \tag6.3.2$$
In fact, we will be using only a lower bound in (6.3.1).

For completeness, let us sketch the main ingredients of the proof of (6.3.1).

Recall that the {\it permanent} of an $N \times N$ matrix $A=\left(a_{ij}\right)$ is defined by
the formula
$$\per A=\sum_{\sigma \in S_N} \prod_{i=1}^N a_{i \sigma(i)},$$
where the sum is taken over all $N!$ permutations $\sigma$ from the symmetric group $S_N$.
For an $m \times n$ matrix $X=\left(x_{ij}\right)$ let us define the $N \times N$ 
block matrix $A(X)$ that has $mn$ blocks of sizes $r_i \times c_j$ for $i=1, \ldots, m$
and $j=1, \ldots, n$ with the $(i,j)$-th block filled by the copies of $x_{ij}$.
A combinatorial computation produces the following expansion
$$\per A(X)=\left(\prod_{i=1}^m  r_i! \right) \left(\prod_{j=1}^n c_j! \right) 
\sum \Sb D=\left(d_{ij}\right) \endSb {x_{ij}^{d_{ij}} \over d_{ij}!},$$
where the sum is taken over all $m \times n$ non-negative integer matrices 
$D=\left(d_{ij}\right)$ with row sums $R$ and column sums $C$.
From this expansion we obtain the formula
$$\#(R,C) =\left(\prod_{i=1}^m {1 \over r_i!} \right) \left(\prod_{j=1}^n {1 \over c_j!} \right) 
\int_{{\Bbb R}^d_+} \per A(X) \exp\left\{-\sum_{ij} x_{ij} \right\} \ dX,$$
where ${\Bbb R}^d_+$ is the set of $m \times n$ positive matrices $X$, see Theorem 1.1 of 
\cite{Ba08}.
Since $\per A(X)$ is a homogeneous polynomial of degree $N$ in $X$, a standard 
change of variables results in the formula 
$$\#(R,C) ={(N+mn-1)! \over \sqrt{mn}}
 \left(\prod_{i=1}^m {1 \over r_i!} \right) \left(\prod_{j=1}^n {1 \over c_j!} \right) 
\int_{\Delta} \per A(X) \ dX,$$
cf. Lemma 4.1 of \cite{Ba08}.
Given a matrix $X \in \Delta$, let $\lambda_1, \ldots, \lambda_m$ and $\mu_1, \ldots, \mu_n$ 
be its scaling factors so that (6.1.1) holds. Let $B(X)$ be the matrix obtained by dividing 
the entries in the $(i,j)$-th block of $A(X)$ by $\lambda_i r_i \mu_j c_j$, so the 
entries in the $(i,j)$-th block of $B(X)$ are equal to $l_{ij}/r_i c_j$.
Hence
$$\per A(X) = \left( \prod_{i=1}^m r_i^{r_i} \right) \left(\prod_{j=1}^n c_j^{c_j} \right) \phi(X) 
\per B(X),$$
cf. Section 3.1 of \cite{Ba08}.
Now we notice that $B(X)$ is a {\it doubly stochastic matrix}, that is, a 
non-negative matrix with row and column sums equal to 1. 
The classical estimate for permanents of doubly stochastic matrices conjectured by van der
Waerden and 
proved by Falikman and Egorychev (see \cite{Fa81}, \cite{Eg81}, and 
Chapter 12 of \cite{LW01}) asserts that
$$\per B(X) \geq {N! \over N^N}$$ 
and hence the lower bound in (6.3.1) follows. The upper bound in (6.3.1) follows from the inequality
for permanents conjectured by Minc and proven by Bregman, 
(see \cite{Br73} and Chapter 11 of \cite{LW01}), which results in
$$\per B(X) \ \leq \ \min\left\{ \prod_{i=1}^m {r_i! \over r_i^{r_i}}, \quad \prod_{j=1}^n {c_j! \over c_j^{c_j}}
\right\},$$
since the entries in the $(i,j)$-th block of $B(X)$ do not exceed
$\min\{1/r_i, \ 1/c_j \}$, see Section 5 of \cite{Ba08} for 
details.

\subhead (6.4) Slicing the simplex \endsubhead
The crucial observation which makes the integral 
$$\int_{\Delta} \phi(X) \ dX$$ 
amenable to analysis is that the simplex $\Delta$ can be sliced by affine subspaces of codimension
$m+n-1$ into sections on which function $\phi$ remains constant.

Let us choose some positive $\lambda_1, \ldots, \lambda_m$ and $\mu_1, \ldots, \mu_n$
and let us consider the affine subspace $L \subset {\Bbb R}^d$ of $m \times n$ matrices
$X=\left(x_{ij}\right)$ satisfying the equations
$$\aligned &\sum_{j=1}^n {x_{ij} \over \lambda_i \mu_j} =r_i \quad \text{for} \quad i=1, \ldots m \\
& \sum_{i=1}^m {x_{ij} \over \lambda_i \mu_j} =c_j \quad \text{for} \quad j=1, \ldots, n
\endaligned \tag6.4.1$$
Clearly, 
$$\phi(X)=\left(\prod_{i=1}^m \lambda_i^{r_i} \right) \left(\prod_{j=1}^n \mu_j^{c_j} \right)
 \quad \text{for all} 
\quad X \in \Delta \cap L. \tag6.4.2$$
Moreover, $\dim L=(m-1)(n-1)$. 

\subhead (6.5) Modification for weighted tables \endsubhead
Similar identities an inequalities hold for weighted tables. For a 
{\it positive} matrix $W=\left(w_{ij}\right)$
of weights, we define the function
$$\phi_{R,C;W}(X)=\phi_{R,C}(Y) \quad \text{where} \quad y_{ij}=w_{ij} x_{ij} \quad 
\text{for all} \quad i,j$$
and $\phi_{R,C}$ is the unweighted function defined in Section 6.2.
Then
$$\split &T(R,C;W) \ \geq \ {N! (N+mn-1)! \over N^N \sqrt{mn}} \\ & \qquad \qquad 
\times \left(\prod_{i=1}^m {r_i^{r_i} \over r_i!}\right)  \left( \prod_{j=1}^n {c_j^{c_j} \over c_j!} \right) 
\int_{\Delta} \phi_{R,C;W}(X) \ dX  \\
& \text{and} \\
&T(R,C:W) \ \leq\  {(N+mn-1)! \over \sqrt{mn}} \\ &\qquad \qquad 
\times \min\left\{ \prod_{i=1}^m {r_i^{r_i} \over r_i!}, \quad
\prod_{j=1}^n {c_j^{c_j} \over c_j!}\right\}  \int_{\Delta} \phi_{R,C;W}(X) \ dX,  \endsplit \tag6.5.1 $$
see \cite{Ba07}, \cite{Ba08}, and the proof sketch in Section 6.3.

Let us choose some positive $\lambda_1, \ldots, \lambda_m$ and $\mu_1, \ldots, \mu_n$
and let us consider the subspace $L \subset {\Bbb R}^d$ of $m \times n$ matrices
$X=\left(x_{ij}\right)$ satisfying the equations
$$\aligned &\sum_{j=1}^n {w_{ij} x_{ij} \over \lambda_i \mu_j} =r_i \quad \text{for} \quad i=1, \ldots m \\
& \sum_{i=1}^m {w_{ij} x_{ij} \over \lambda_i \mu_j} =c_j \quad \text{for} \quad j=1, \ldots, n
\endaligned \tag6.5.2$$
Clearly, 
$$\phi_{R,C;W}(X)=\left(\prod_{i=1}^m \lambda_i^{r_i} \right) \left(\prod_{j=1}^n \mu_j^{c_j} \right)
 \quad \text{for all} 
\quad X \in \Delta \cap L. \tag6.5.3$$
Moreover, $\dim L=(m-1)(n-1)$.

\head 7. Proofs of Theorems 1.1 and 1.3 and Lemma 1.4 \endhead

We prove Lemma 1.4 first.
\subhead (7.1) Proof of Lemma 1.4 \endsubhead It is straightforward to check that the function
$$g(x; w)=(x+1) \ln (x+1) - x\ln x +x \ln w \quad \text{for} \quad x \geq 0$$
is strictly concave for $x>0$. Therefore, the maximum of $g(X;W)$ on $\PP(R,C)$ is attained 
at a single point $Z=\left(z_{ij}\right)$. Let us show that necessarily $z_{ij}>0$ for all $i,j$.

Since 
$$ g'(x; w)= \ln \left({x+1 \over x}\right) + \ln w,$$
the derivative of $g(x; w)$ at $x>0$ is finite and the right derivative at  $x=0$ is $+\infty$.
Let $Y \in \PP(R,C)$ be a matrix with positive entries,
for example, $Y=\left(y_{ij}\right)$ where $y_{ij} =r_i c_j/N$. If $z_{ij}=0$ for some $i,j$ then
$$g\bigl((1-\epsilon) Z +\epsilon Y; W \bigr) \ > \ g\bigl(Z; W \bigr)$$
for some sufficiently small $\epsilon >0$, which is a contradiction.

Thus $z_{ij}>0$ for all $i,j$ and hence $Z$ lies in the relative interior of $\PP(R,C)$.
Therefore the gradient of $g(X; W)$ at $X=Z$ is orthogonal to the affine span of $\PP(R,C)$,
that is,
$$\ln \left({z_{ij}+1 \over z_{ij}}\right) + \ln w_{ij} =\lambda_i +\mu_j \quad \text{for all} \quad i,j$$
and some $\lambda_1, \ldots, \lambda_m$, $\mu_1, \ldots, \mu_n$.

Let 
$$\xi_i=e^{-\lambda_i} >0 \quad \text{for} \quad i=1, \ldots, m \quad \text{and}
\quad \eta_j=e^{-\mu_j} >0 \quad \text{for} \quad j=1, \ldots, n.$$
Then 
$$w_{ij} \xi_i \eta_j ={z_{ij} \over z_{ij}+1} < 1 \quad \text{for all} \quad i,j,$$
and
$$z_{ij}={w_{ij} \xi_i \eta_j \over 1-w_{ij} \xi_i \eta_j} \quad \text{for all} \quad i,j. \tag7.1.1$$
In particular, 

$$\split &\sum_{j=1}^n {w_{ij} \xi_i \eta_j \over 1-w_{ij} \xi_i \eta_j} =r_i \quad 
\text{for} \quad i=1, \ldots, m \quad \text{and} \\
&\sum_{i=1}^m {w_{ij} \xi_i \eta_j \over 1-w_{ij} \xi_i \eta_j }=c_j \quad 
\text{for} \quad j=1, \ldots, n. \endsplit \tag7.1.2$$

Equations (7.1.2) are equivalent to the statement that the point \break
$\ttt^{\ast}=\left(\lambda_1, \ldots, \lambda_m\right)$ and $\sss^{\ast}=\left(\mu_1, \ldots, \mu_n \right)$
is a critical point of the function
$$\phi(\ttt, \sss) =\sum_{i=1}^m r_i t_i + \sum_{j=1}^n c_j s_j 
-\sum_{ij} \ln \left(1-w_{ij} e^{-t_i -s_j} \right).$$
Since $\phi$ is convex, the point $\left(\sss^{\ast}, \ttt^{\ast} \right)$ is a minimum point of 
$\phi$ and hence the point $\xx^{\ast}=\left(\xi_1, \ldots, \xi_m\right)$ and 
$\yy^{\ast}=\left(\eta_1, \ldots, \eta_n \right)$ is a point where the infimum of 
$$F(\xx, \yy; W) =\left(\prod_{i=1}^m x_i^{-r_i} \right)  \left(\prod_{j=1}^n y_j^{-c_j} \right) 
\left(\prod_{ij} {1 \over 1-w_{ij} x_i y_j}\right)$$
is attained in the region $x_1, \ldots, x_m >0$, $y_1, \ldots, y_n>0$, and $w_{ij} x_i y_j <1$
for all $i,j$. 

Using (7.1.1) and (7.1.2), we conclude that 
$$\split g(Z; W) =&\sum_{ij} (z_{ij}+1) \ln (z_{ij}+1) -\sum_{ij} z_{ij}(\ln z_{ij} - \ln w_{ij})\\=
&-\sum_{ij} {\ln (1-w_{ij} \xi_i \eta_j) \over 1-w_{ij} \xi_i \eta_j} -\sum_{ij} 
{w_{ij} \xi_i \eta_j \over 1-w_{ij} \xi_i \eta_j} \ln \left({ \xi_i \eta_j \over 1-w_{ij} \xi_i \eta_j}\right)\\
=& -\sum_{ij} \ln (1-w_{ij} \xi_i \eta_j ) -\sum_{i=1}^m \ln \xi_i 
\left(\sum_{j=1}^n {w_{ij} \xi_i \eta_j \over 1-w_{ij} \xi_i \eta_j} \right) 
\\& \qquad \qquad \qquad -\sum_{j=1}^n \ln \eta_j \left(\sum_{j=1}^n {w_{ij} \xi_i \eta_j \over 1-w_{ij} \xi_i \eta_j} \right) \\=
&-\sum_{i=1}^m r_i \xi_i -\sum_{j=1}^n c_j \eta_j - \sum_{ij} \ln(1-w_{ij} \xi_i \eta_j)=
\ln F(\xx^{\ast}, \yy^{\ast}; W),
\endsplit$$
as claimed.

We observe that the value of $F(\xx, \yy; W)$ does not change if we scale 
$x_i \longmapsto x_i \tau$, $y_j \longmapsto y_j \tau^{-1}$ for $\tau > 0$. In the case of 
$w_{ij}=1$ for all $i,j$ we have $\xi_i \eta_j <1$ for all $i,j$ and hence by choosing an 
appropriate $\tau$ we can enforce $0< \xi_i, \eta_j <1$ for all $i,j$.
{\hfill \hfill \hfill} \qed
\bigskip

We consider the space ${\Bbb R}^d$ for $d=mn$ of $m \times n$ real matrices, the 
affine hyperplane $\AA \subset {\Bbb R}^d$ defined by the equation 
$\sum_{ij} x_{ij}=1$ and the standard open simplex $\Delta \subset \AA$ defined by the 
inequalities $x_{ij} > 0$ for all $i,j$. Let $\phi=\phi_{R,C;W}$ be the function defined in
Sections 6.5 and 6.2.

We start with a technical lemma, which is a straightforward modification of Lemma 4.3.
\proclaim{(7.2) Lemma} Let $L \subset \AA$ be an affine subspace, $\dim L=d-k-1$ for  $k\geq 1$.
Suppose that there is 
a point $A=\left(\alpha_{ij}\right)$, $A \in \Delta \cap L$, such that 
$$\alpha_{ij} \geq {\epsilon \over d} \quad \text{for all} \quad i,j.$$
Suppose further that the value of the function $\phi=\phi_{R,C;W}$ on $\Delta \cap L$ 
is constant and equal to $\tau$.
Then 
$$\int_{\Delta} \phi(X) \ dX \geq \gamma \left({2 \epsilon \over d^2 (N+1)} \right)^k 
\tau \vl_{d-k-1}(\Delta \cap L)$$
for some absolute constant $\gamma>0$ (one can choose $\gamma =e^{-2} \approx 0.14$).
\endproclaim
\demo{Proof}
Let 
$$Q_0=\left\{ {1 \over d}A +{d-1 \over d}X: \quad X \in \Delta \cap L \right\}.$$
As in the proof of Lemma 4.3, we have 
$$\vl_{d-k-1} (Q_0) \geq {1 \over e} \vl_{d-k-1} (\Delta \cap L)$$ 
and for any $X \in Q_0$, $X =\left(x_{ij}\right)$, we have 
$$x_{ij} \geq {\epsilon \over d^2} \quad \text{for all} \quad i,j.$$
Let us define $Q$ by
$$Q=\left\{X \in \Delta: \quad \|X-Y\|_{\infty} \leq {\epsilon \over d^2 (N+1)}  \quad \text{for some} \quad 
Y \in Q_0 \right\}.$$
Then, as in Lemma 4.3, we have 
$$\vl Q \ \geq \ {1 \over e} \left({2\epsilon \over d^2 (N+1)} \right)^k \vl_{d-k-1}(\Delta \cap L).$$
We note that for every $X \in Q$ there is a $Y \in \Delta \cap L$ such that 
$$x_{ij} \geq \left(1- {1 \over N+1} \right) y_{ij} \quad \text{for all} \quad i,j.$$
Since $\phi$ is monotone and homogeneous of degree $N$ (see Section 6.2) , 
we have 
$$\phi(X) \ \geq \ \left(1- {1 \over N+1} \right)^N \tau \ > \ {1 \over e} \tau \quad \text{for all} \quad X \in Q.$$
Since 
$$\int_{\Delta} \phi(X) \ dX \geq \int_Q \phi(X) \ d X,$$ 
the proof follows.
{\hfill \hfill \hfill} \qed
\enddemo

\subhead (7.3) Proof of Theorem 1.1 \endsubhead The upper bound follows immediately from the 
standard generating function expression:
$$\prod_{i,j} {1 \over 1-x_i y_j} =\sum_{R,C} \#(R,C)x^R y^C, 
\quad \text{where} \quad x^R =\prod_{i=1}^m x_i^{r_i} \quad \text{and} \quad 
x^C =\prod_{j=1}^n y_j^{c_j}$$ 
and the some is taken over all pairs of positive integer $m$-vectors $R=\left(r_1, \ldots, r_m\right)$ and $n$-vectors $C=\left(c_1, \ldots, c_n \right)$
such that $ r_1+ \ldots +r_m =c_1 + \ldots + c_n$. 

Let us prove the lower bound. By Lemma 1.4 the minimum of  
$$F(\xx,\yy)=\left(\prod_{i=1}^m x_i^{-r_i}\right)\left( \prod_{j=1}^n y_j^{-c_j} \right)
\left(\prod_{ij} {1 \over 1-x_i y_j}\right)$$
on the open cube $0<x_i, y_j <1$ for all $i,j$ is  attained at a certain point
$$\xx^{\ast}=\left(\xi_1, \ldots, \xi_m\right) \quad \text{and} \quad 
\yy^{\ast}=\left(\eta_1, \ldots, \eta_n \right),$$
which, moreover, satisfies 
$$\aligned &\sum_{j=1}^n {\xi_i \eta_j \over 1-\xi_i \eta_j} =r_i \quad \text{for} \quad i=1, \ldots, m \\
&\sum_{i=1}^m {\xi_i \eta_j \over 1-\xi_i \eta_j} =c_j \quad \text{for} \quad j=1, \ldots, n. \endaligned
\tag7.3.1$$
Equations (7.3.1) can also be obtained by 
setting the gradient of $\ln F$ to 0.

In the space of $m \times n$ matrices ${\Bbb R}^d$ with $d=mn$, let us consider the standard 
simplex $\Delta$ and the point $A=\left(\alpha_{ij}\right)$ defined by 
$$\alpha_{ij}= {1 \over (N +mn) (1-\xi_i \eta_j)} \quad \text{for all} \quad i,j.$$
By (7.3.1), we have 
$$\sum_{j=1}^n {1 \over 1-\xi_i \eta_j} = \sum_{j=1}^n {1 -\xi_i \eta_j \over 1-\xi_i \eta_j} +
\sum_{j=1}^n {\xi_i \eta_j \over 1-\xi_i \eta_j} = n + r_i \quad \text{for all} \quad i, $$
so $A$ lies in $\Delta$. 
Let 
$$\split &\lambda_i = {1 \over \xi_i \sqrt{N+mn}} \quad \text{for} \quad i=1, \ldots, m \quad \text{and}\\
&\mu_j={1 \over \eta_j  \sqrt{N+mn}} \quad \text{for} \quad j=1, \ldots, n. \endsplit$$
Let us consider the affine subspace $L \subset {\Bbb R}^d$ defined by the 
system of equations
$$\split &\sum_{j=1}^n {x_{ij} \over \lambda_i \mu_j} =r_i \quad \text{for} \quad i=1, \ldots m \\
&\sum_{i=1}^m {x_{ij} \over \lambda_i \mu_j} =c_j \quad \text{for} \quad j=1, \ldots, n. \endsplit$$
Hence $\dim L =(m-1)(n-1)$ and $A \in L$ by (7.3.1).

By (6.4.2), the density $\phi=\phi_{R,C}$ is constant on $L$ and equal to 
$$\tau = {1 \over (N+mn)^N} \left(\prod_{i=1}^m \xi_i^{-r_i}\right)\left( \prod_{j=1}^n \eta_j^{-c_j}\right).$$
By Part (1) of Theorem 3.1, the volume of the section $\Delta \cap L$ within a factor 
of $(N+mn)^{O(m+n)}$
is at least 
$$ {e^{mn} \over (N+mn)^{mn}} \left(\prod_{ij} {1 \over 1-\xi_i \eta_j}\right).$$
More precisely, for $k=d-1-\dim (\Delta \cap L)$ we have $k=m+n-1$ or $k=m+n-2$ and
$$\vl_{d-1-k} (\Delta \cap L) \ \geq \ {\Gamma(k/2+1) \over 2e^3  \pi^{k/2} \sqrt{mn}} 
{(mn)^{mn} \over (mn)! }{1 \over (N+mn)^{mn}} \left(\prod_{ij} {1 \over 1-\xi_i \eta_j}\right).$$

Choosing $\epsilon=1/(N+mn)$ in Lemma 7.2, we estimate the integral
$$\int_{\Delta} \phi(X) \ dX$$ within a factor of $N^{O(m+n)}$ from below by 
$${e^{mn} \over (N+mn)^{N+mn}} \left(\prod_{i=1}^m \xi_i^{-r_i} \right)
\left( \prod_{j=1}^n \eta_j^{-c_j} \right)
\left(\prod_{ij} {1 \over 1-\xi_i \eta_j}\right).$$ 
More precisely,
$$\split \int_{\Delta} \phi(X) \ dX \ \geq \ &{\Gamma\left({m+n \over 2} \right) \over 2e^5  
\pi^{{m+n-2 \over 2}} \sqrt{mn}} \left({2\over (mn)^2 (N+1) (N+mn)}\right)^{m+n-1} \\
& \qquad \times {(mn)^{mn} \over (mn)! }{1 \over (N+mn)^{N+mn}} \\
&\qquad \times \left(\prod_{i=1}^m \xi_i^{-r_i} \right) \left(\prod_{j=1}^n \eta_j^{-c_j}\right)
 \left(\prod_{ij} {1 \over 1-\xi_i \eta_j}\right)\endsplit$$
provided $m + n \geq 10$.
Hence by (6.3.2) the number $\#(R, C)$ is estimated from below 
within a factor of $(N+mn)^{O(m+n)}$ by 
$$\split e^N (N+mn)!  \int_{\Delta} \phi(X) \ dX \ \approx \ &{e^{mn+N} (N+mn)!\over (N+mn)^{N+mn}} 
\\
& \qquad \times \left(\prod_{i=1}^m \xi_i^{-r_i} \right)\left(\prod_{j=1}^n \eta_j^{-c_j} \right)
\left(\prod_{ij} {1 \over 1-\xi_i \eta_j} \right) \\ \approx &  F(x^{\ast}, y^{\ast})=\rho,\endsplit$$
where ``$\approx$'' stands for an approximation within a $N^{O(m+n)}$ factor.

More precisely, by (6.3.1)
$$\split \#(R,C) \ \geq \  \ &{\Gamma\left({m+n \over 2} \right) \over 2e^5  
\pi^{{m+n-2 \over 2}} mn (N+mn)} \left({2\over (mn)^2 (N+1) (N+mn)}\right)^{m+n-1} \\
&\quad \times \left(\prod_{i=1}^m {r_i^{r_i} \over r_i!} \right) \left( \prod_{j=1}^n {c_j^{c_j} \over c_j!} \right) 
{N! (N+mn)! (mn)^{mn} \over  N^N  (N+mn)^{N+mn}(mn)!} \rho(R,C) \endsplit$$
provided $m +n \geq 10$.
{\hfill \hfill \hfill} \qed

The proof of Theorem 1.3 is a straightforward modification of the proof of Theorem 1.1.

\subhead (7.4) Proof of Theorem 1.3 \endsubhead
The upper bound follows from the generating function 
expression 
$$\prod_{ij} {1 \over 1-w_{ij} x_i y_j} =\sum_{R, C} T(R, C; W) x^R y^C.$$
Let us prove the lower bound. Since $T(R,C;W)$ is a polynomial in $W$, without loss of generality we assume that $W$ is a strictly 
positive matrix. Let 
$$x^{\ast}
=\left(\xi_1, \ldots, \xi_m \right) \quad \text{and} \quad y^{\ast}=\left(\eta_1, \ldots, \eta_n \right)$$ 
be the minimum point of 
$$F(\xx,\yy; W)=\left(\prod_{i=1}^m x_i^{-r_i} \right) \left(\prod_{j=1}^n y_j^{-c_j}\right)
\left( \prod_{ij} {1 \over 1-w_{ij} x_i y_j}\right),$$ 
see Lemma 1.4.
Then
$$\aligned 
&\sum_{j=1}^n {w_{ij} \xi_i \eta_j \over 1-w_{ij} \xi_i \eta_j} =r_i \quad \text{for} \quad i=1, \ldots, m
\\ &\sum_{i=1}^m {w_{ij} \xi_i \eta_j \over 1-w_{ij}\xi_i \eta_j} =c_j \quad \text{for} \quad j=1, \ldots, n.
\endaligned \tag7.4.1$$
In the space of matrices, let us consider the standard simplex $\Delta$ and the matrix 
$A=\left(\alpha_{ij}\right)$ 
$$\alpha_{ij}={1 \over (N+mn)(1-w_{ij}\xi_i \eta_j)} \quad \text{for all} \quad i,j.$$
As in the proof of Theorem 1.1, we check from (7.4.1) that indeed $A \in \Delta$. 
Let 
$$\split &\lambda_i={1 \over \xi_i \sqrt{N +mn}} \quad \text{for} \quad i=1, \ldots, m \quad \text{and} \\
&\mu_j={1 \over \eta_j \sqrt{N+mn}} \quad \text{for} \quad j=1, \ldots, n. \endsplit$$ 

Let us consider the affine space $L \subset {\Bbb R}^d$ defined by the equations 
$$\split & \sum_{j=1}^n {w_{ij} x_{ij} \over \lambda_i \mu_j} =r_i \quad \text{for} \quad i=1, \ldots, m \\
&\sum_{i=1}^m {w_{ij} x_{ij} \over \lambda_i \mu_j} =c_j \quad \text{for} \quad j=1, \ldots, n .
\endsplit$$
Then $A \in L$, the value of
$\phi_{R,C;W}$ on $\Delta \cap L$ is constant and equal to
$$\tau={1 \over (N+mn)^N} \left(\prod_{i=1}^m \xi_i^{-r_i} \right) \left(\prod_{j=1}^n \eta_j^{-c_j}\right),$$
 see (6.5.2)-(6.5.3). Next, we use the lower bound in (6.5.1) and the proof proceeds as for Theorem 1.1.
{\hfill \hfill \hfill} \qed

\head Acknowledgments \endhead

This work grew out of a joint project with Alex Samorodnitsky, and at an earlier stage also with 
Alexander Yong,  on constructing computationally efficient algorithms for enumeration of contingency tables, see \cite{BSY07}. I am grateful to Keith Ball and Matthieu Fradelizi for teaching
me methods to estimate volumes of sections of convex bodies. I benefitted from 
conversations with Alex Samorodnitsky, Imre B\'ar\'any, and Roy Meshulam.


\Refs
\widestnumber\key{AAAAA}

\ref\key{BLV04}
\by W. Baldoni-Silva, J.A. De Loera, and M. Vergne
\paper Counting integer flows in networks
\jour Found. Comput. Math. 
\vol 4 
\yr 2004
\pages 277--314
\endref


\ref\key{Bal88}
\by K. Ball
\paper Logarithmically concave functions and sections of convex sets in $ R\sp n$
\jour Studia Math. 
\vol 88 
\yr 1988
\pages 69--84
\endref

\ref\key{Bal97}
\by K. Ball
\paper An elementary introduction to modern convex geometry
\inbook  Flavors of Geometry
\pages 1--58
\bookinfo Math. Sci. Res. Inst. Publ., 31
\publ Cambridge Univ. Press
\publaddr Cambridge
\yr 1997
\endref

\ref\key{Ba07}
\by A. Barvinok
\paper  Brunn-Minkowski inequalities for contingency tables and integer flows
\jour Adv. Math. 
\vol 211 
\yr 2007 
\pages 105--122
\endref

\ref\key{Ba08}
\by A. Barvinok
\paper Enumerating contingency tables via random permanents \newline
\jour Combin. Probab. Comput. 
\vol 17 
\yr 2008
\pages 1--19
\endref

\ref\key{BSY07}
\by A. Barvinok, A, Samorodnitsky, and A. Yong
\paper Counting magic squares in quasi-polynomial time
\paperinfo preprint, arXiv math.CO/0703227 
\yr 2007
\endref

\ref\key{BP03}
\by M. Beck and D. Pixton
\paper The Ehrhart polynomial of the Birkhoff polytope
\jour Discrete Comput. Geom. 
\vol 30 
\yr 2003
\pages  623--637
\endref 

\ref\key{BBK72}
\by A. B\' ek\' essy, P.  B\'  ek\' essy, and  J. Koml\' os
\paper Asymptotic enumeration of regular matrices 
\jour Studia Sci. Math. Hungar. 
\vol 7 
\yr 1972
\pages  343--353
\endref 

\ref\key{Ben74}
\by E. Bender
\paper The asymptotic number of non-negative integer matrices with given row and column sums 
\jour Discrete Math. 
\vol 10 
\yr 1974
\pages 217--223
\endref 

\ref\key{Bol97}
\by B. Bollob\' as
\paper Volume estimates and rapid mixing
\inbook  Flavors of Geometry
\pages 151--182
\bookinfo Math. Sci. Res. Inst. Publ., 31
\publ Cambridge Univ. Press
\publaddr Cambridge
\yr 1997
\endref 

\ref\key{Br73}
\by L.M. Bregman 
\paper Certain properties of
nonnegative matrices and their permanents 
\jour Dokl. Akad. Nauk
SSSR 
\vol 211 
\yr 1973
\pages 27--30 
\endref

\ref\key{CM07a}
\by E.R. Canfield and B.D. McKay
\paper Asymptotic enumeration of contingency tables with constant margins
\paperinfo  preprint arXiv math.CO/0703600 
\yr 2007
\endref

\ref\key{CM07b}
\by E.R. Canfield and B.D. McKay
\paper The asymptotic volume of the Birkhoff polytope
\paperinfo preprint arXiv:0705.2422 
\yr 2007
\endref

\ref\key{C+05}
\by Y. Chen, P. Diaconis, S.P.  Holmes, and J.S. Liu
\paper Sequential Monte Carlo methods for statistical analysis of tables
\jour J. Amer. Statist. Assoc. 
\vol 100 
\yr 2005
\pages 109--120
\endref

\ref\key{CD03}
\by M. Cryan and M. Dyer
\paper A polynomial-time algorithm to approximately count contingency tables when the number of rows is constant
\paperinfo  Special issue on STOC2002 (Montreal, QC)
\jour J. Comput. System Sci.
\vol 67 
\yr 2003
\pages  291--310
\endref

\ref\key{DE85}
\by P. Diaconis and B. Efron
\paper Testing for independence in a two-way table: new interpretations of the chi-square statistic. With discussions and with a reply by the authors
\jour Ann. Statist. 
\vol 13 
\yr 1985
\pages 845--913
\endref

 \ref\key{DG04}
 \by P. Diaconis and A. Gamburd
 \paper Random matrices, magic squares and matching polynomials
 \jour Electron. J. Combin. 
 \vol 11 
 \yr 2004/06
 \paperinfo Research Paper 2, 26 pp. (electronic)
 \endref

\ref\key{DG95}
\by P. Diaconis and A. Gangolli
\paper Rectangular arrays with fixed margins
\inbook Discrete probability and algorithms (Minneapolis, MN, 1993)
\pages 15--41 
\bookinfo IMA Vol. Math. Appl., 72
\publ  Springer
\publaddr New York
\yr 1995
\endref

\ref\key{DLY07}
\by J.A. De Loera, F. Liu, and R. Yoshida
\paper Formulas for the volumes of the polytope of doubly-stochastic matrices and its faces
\paperinfo preprint arXiv math.CO/0701866 
\yr 2007
\endref

 \ref\key{DKM97}
 \by M. Dyer, R. Kannan, and J. Mount
 \paper Sampling contingency tables
 \jour Random Structures Algorithms 
 \vol 10 
 \yr 1997
 \pages 487--506
 \endref

\ref\key{Eg81}
\by G.P. Egorychev 
\paper The solution of van der
Waerden's problem for permanents 
\jour Adv. in Math. 
\vol 42
\yr 1981
\pages 299--305
\endref

\ref\key{Fa81} 
\by D.I. Falikman 
\paper Proof of the van der
Waerden conjecture on the permanent of a doubly stochastic matrix
(Russian) 
\jour Mat. Zametki 
\vol 29 
\yr 1981 
\pages 931--938
\endref

\ref\key{Frad97}
\by M. Fradelizi
\paper Sections of convex bodies through their centroid
\jour Arch. Math. (Basel)
\vol 69 
\yr 1997
\pages 515--522
\endref 

\ref\key{Goo76}
\by I.J. Good
\paper On the application of symmetric Dirichlet distributions and their mixtures to contingency tables
\jour Ann. Statist. 
\vol 4 
\yr 1976
\pages 1159--1189
\endref

\ref\key{GM07}
\by C. Greenhill and B.D. McKay
\paper  Asymptotic enumeration of sparse nonnegative integer matrices with specified row and column sums
\paperinfo preprint arXiv:0707.0340 
\yr  2007
\endref

\ref\key{Gr\"u60}
\by B. Gr\"unbaum
\paper Partitions of mass-distributions and of convex bodies by hyperplanes
\jour Pacific J. Math. 
\vol10 
\yr 1960 
\pages 1257--1261
\endref

\ref\key{Khi57}
\by A.I. Khinchin
\book Mathematical Foundations of Information Theory
\bookinfo Translated by R. A. Silverman and M. D. Friedman
\publ  Dover Publications, Inc.
\publaddr New York, N. Y.
\yr 1957
\endref

\ref\key{LW01}
 \by J.H. van Lint and R.M. Wilson 
 \book A Course in
Combinatorics. Second edition 
\publ Cambridge University Press
\publaddr Cambridge 
\yr 2001 
\endref

\ref\key{MO68}
\by A.W. Marshall and I. Olkin
\paper Scaling of matrices to achieve specified row and column sums
\jour Numer. Math. 
\vol 12 
\yr 1968 
\pages 83--90
\endref 

\ref\key{Mor02}
\by B.J. Morris
\paper  Improved bounds for sampling contingency tables
\jour Random Structures Algorithms
\vol 21 
\yr 2002
\pages 135--146
\endref

\ref\key{NN94}
\by Yu. Nesterov and A. Nemirovskii
\book Interior-Point Polynomial Algorithms in Convex Programming
\bookinfo SIAM Studies in Applied Mathematics, 13
\publ Society for Industrial and Applied Mathematics (SIAM)
\publaddr Philadelphia, PA
\yr 1994
\endref

\ref\key{Pak00}
\by I. Pak
\paper Four questions on Birkhoff polytope
\jour Ann. Comb. 
\vol 4 
\yr 2000
\pages  83--90
\endref

\ref\key{RS89}
\by U.G. Rothblum and H. Schneider
\paper Scalings of matrices which have prespecified row sums and column sums via optimization
\jour Linear Algebra Appl. 
\vol 114/115 
\yr 1989
\pages  737--764
\endref

\ref\key{Sch92}
\by M. Schmuckenschl\" ager
\paper On the volume of the double stochastic matrices 
\jour Acta Math. Univ. Comenian. (N.S.) 
\vol 61 
\yr 1992
\pages 189--192
\endref 

\ref\key{Si64}
\by R. Sinkhorn
\paper A relationship between arbitrary positive matrices and doubly stochastic matrices
\jour Ann. Math. Statist. 
\vol 35 
\yr 1964 
\pages 876--879
\endref

\ref\key{Vaa79}
\by J.D. Vaaler
\paper A geometric inequality with applications to linear forms
\jour Pacific J. Math.
\vol 83 
\yr 1979
\pages 543--553
\endref 

\endRefs
\enddocument
\end